\documentclass{amsart}

\usepackage{amsmath,amssymb}
\usepackage[arrow,matrix]{xy}

\theoremstyle{proclaim}

\newtheorem*{Th*}{Theorem}

\newtheorem*{Cor*}{Corollary}

\theoremstyle{statement}

\theoremstyle{fancyproclaim}

\numberwithin{equation}{section}

\makeatletter
\def\Set@Scallop[#1]#2#3{{#1}\Parens{#2}{#3}}
\newcommand\DeclareScalableOperator[2]{%
  \expandafter\def\csname#1\endcsname{\@ifnextchar[{{#2}\Set@Scallop}{{#2}\Set@Scallop[{}]}}
}
\makeatother

\DeclareScalableOperator{Ct}{\mathcal{C}} 
\DeclareScalableOperator{Cc}{\mathcal{K}} 
\DeclareScalableOperator{Lp}{\mathbf{L}} 
\DeclareScalableOperator{Hom}{\mathrm{Hom}} 
\DeclareScalableOperator{Todd}{\mathrm{Td}} 

\newcommand{\fa}{for all }
\newcommand\mathfa[1][{}]{\quad\text{\fa{#1} }}
\newcommand\smathfa[1][{}]{\ \text{\fa{#1} }}

\newcommand{\scth}{such that }
\newcommand{\AND}{and}

\newcommand\mathtxt[1]{\quad\text{{#1}}\quad}

\newcommand{\nd}{\mathtxt\AND}

\newcommand\vphi{\varphi}

\newcommand\eps{\varepsilon}

\newcommand\ints{\mathbb{Z}}
\newcommand\rats{\mathbb{Q}}
\newcommand\reals{\mathbb{R}}
\newcommand\cplxs{\mathbb{C}}
\newcommand\knums{\mathbb K}
\newcommand\sph{\mathbb{S}}
\newcommand\ball{\mathbb{B}}

\newcommand\sle{\leqslant}
\newcommand\sge{\geqslant}

\DeclareMathOperator\Index{\mathrm{index}}
\DeclareMathOperator\ch{\mathrm{ch}}
\DeclareMathOperator\rk{\mathrm{rk}}

\makeatletter
\newcommand\Size[7][1]{
                                 \ifx#20%
                                        \def\r@l{}\def\r@m{}\def\r@r{}%
                                 \else%
                                    \ifx#21%
                                           \def\r@l{\bigl}\def\r@r{\bigr}\def\r@m{\bigm}%
                                    \else%
                                           \ifx#22%
                                                 \def\r@l{\Bigl}\def\r@r{\Bigr}\def\r@m{\Bigm}%
                                            \else%
                                                 \ifx#23%
                                                        \def\r@l{\biggl}\def\r@r{\biggr}\def\r@m{\biggm}%
                                                  \else
                                                        \ifx#24%
                                                              \def\r@l{\Biggl}\def\r@r{\Biggr}\def\r@m{\Biggm}%
                                                        \fi%
                                                  \fi%
                                            \fi%
                                      \fi%
                                 \fi%
                                 \ifx#10%
                                       \def\r@m{}%
                                 \fi%
                                 \r@l#3{#4}\r@m#5{#6}\r@r#7%
}%
\makeatother

\newcommand\Set[3]{
                                 \Size{#1}{\{}{#2}{|}{#3}{\}}%
}%
\newcommand\Scp[3]{
                                 \Size{#1}{(}{#2}{|}{#3}{)}%
}%
\newcommand\Rscp[3]{
                                 \Size[0]{#1}{(}{#2}{:}{#3}{)}%
}%
\newcommand\Parens[2]{
  \Size[0]{#1}{(}{#2}{}{}{)}
}

\newcommand\Norm[2]{
  \Size[0]{#1}{\lVert}{#2}{}{}{\rVert}
}
\newcommand\Abs[2]{
  \Size[0]{#1}{\lvert}{#2}{}{}{\rvert}
}
\newcommand\Span[2]{
  \Size[0]{#1}{\langle}{#2}{}{}{\rangle}
}

\makeatletter

\newcommand{\IfUpperCase}[1]{\begingroup 
  \protected@edef\@tempa{\expandafter\@firstofone\@firstofone#1.}%
  \expandafter\IfUpperCasE \@tempa\delimiter}

\def\IfUpperCasE #1#2\delimiter{%
  \protected@edef\@tempa{\meaning#1\meaning a}%
  \ifnum \expandafter\IfUppercaSE\@tempa \IfUppercaSE
   \endgroup \expandafter\@firstoftwo
  \else
   \endgroup \expandafter\@secondoftwo
  \fi}

\@ifundefined{strip@prefix}{\def\strip@prefix#1>{}}{}
\def\@tempa{the letter }
\edef\@tempa{\expandafter\strip@prefix\meaning\@tempa}

\expandafter\def\expandafter\IfUppercaSE\expandafter#\expandafter1\@tempa#2#3\IfUppercaSE{\uccode`#2=`#2 }

\newif\ifuc@se
\def\setuc@se#1{\IfUpperCase{#1}{\uc@setrue}{\uc@sefalse}}

\def\theoremn@me#1{\ifuc@se \lowercase{\csname#1name\endcsname}\ignorespaces%
  \else \edef\@temp{\lowercase{\lowercase{\csname#1name\endcsname}}}\@temp\ignorespaces%
  \fi}
\def\theoremn@mes#1{\ifuc@se \lowercase{\csname#1names\endcsname}\ignorespaces%
  \else \edef\@temp{\lowercase{\lowercase{\csname#1names\endcsname}}}\@temp\ignorespaces%
  \fi}
  
\def\thmref#1#2{\setuc@se{#1}\lowercase{{\theoremn@me{#1}\lowercase{\ref{#1:#2}}}}}

\newcommand{\DefTheorem}[2]{\newenvironment{#1}[1][\empty]{\ignorespaces\begin{#2}\ifx##1\empty{}\else\lowercase{\label{#1:##1}}\fi\ignorespaces}{\end{#2}\ignorespacesafterend}}

\DefTheorem{Th}{theorem}
\DefTheorem{Prop}{proposition}
\DefTheorem{Cor}{corollary}
\DefTheorem{Lem}{lemma}
\DefTheorem{Def}{definition}
\DefTheorem{Rem}{remark}
\DefTheorem{Par}{para}

\newenvironment{Par*}{\ignorespaces\noindent\ignorespaces}{\ignorespacesafterend}
                         
\makeatother

\begin{document}

\date{}

\title[Polytopes and Index]{Convex Polytopes and\\ the Index of Wiener--Hopf Operators}

\author[A.~Alldridge]{Alexander Alldridge}

\address{A.~ALLDRIDGE, Institut f\"ur Mathematik, Universit\"at Paderborn, 33098 Paderborn, Germany}
\email{alldridg@math.upb.de}

\begin{abstract} 
We study the C$^*$-algebra of Wiener--Hopf operators $A_\Omega$ on a cone $\Omega$ with polyhedral base $P\,$. As is known, a sequence of symbol maps may be defined, and their kernels give a filtration by ideals of $A_\Omega\,$, with liminary subquotients. One may define $K$-group valued `index maps' between the subquotients. These form the $E^1$ term of the Atiyah--Hirzebruch type spectral sequence induced by the filtration. We show that this $E^1$ term may, as a complex, be identified with the cellular complex of $P\,$, considered as CW complex by taking convex faces as cells. It follows that $A_\Omega$ is $KK$-contractible, and that $A_\Omega/\knums$ and $S$ are $KK$-equivalent. Moreover, the isomorphism class of $A_\Omega$ is a complete invariant for the combinatorial type of $P\,$. 
\end{abstract}

\maketitle

\section*{INTRODUCTION}

\noindent
The classical \emph{Wiener--Hopf algebra} $A_{\reals_+}\,$, also known as the \emph{reduced Toeplitz algebra}, is an object of basic interest in Operator Theory. It may be defined as the C$^*$-algebra of bounded operators on $\Lp[^2]0{\reals_+}$ generated by all \emph{Wiener--Hopf operators}
\[
(W_fg)(x)=\int_0^\infty f(x-y)g(y)\,dy\mathfa g\in\Lp[^2]0{\reals_+}\,,\,f\in\Lp[^1]0{\reals}\,,\,x\in\reals_+\ .
\]

The \emph{symbol map} $\sigma:A_{\reals_+}\to\Ct[_0]0\reals$ is the surjective $*$-morphism defined on generators by $\sigma(W_f)=\Hat f\,$, where the latter denotes the Fourier transform of $f\in\Lp[^1]0\reals\,$. It gives rise to a short exact sequence
\[
\xymatrix{0\ar[r]&\knums\ar[r]&A_{\reals_+}\ar[r]^-{\sigma}&S=\Ct[_0]0\reals\ar[r]&0}\ .\tag{$*$}
\]

Let $\partial:K_1(S)\to K_0(\knums)=\ints$ be the connecting map in $K$-theory induced by this exact sequence. The following theorem is well-known.

\begin{Th*}
	If $T$ is a Fredholm operator contained in the unitisation of $A_{\reals_+}\,$, then 
	\[
		\partial[\sigma(T)]=\Index T\in\ints=K_0(\knums)\ .
	\]
	Moreover, $\partial$ is a group isomorphism, and for any $n\in\ints\,$, there exists a Fredholm element $T$ as above \scth $\Index T=n\,$. 
\end{Th*}

\noindent
The six-term exact sequence immediately gives the following corollary.

\begin{Cor*}
	We have $K_i(A_{\reals_+})=0$ for $i=0,1\,$. 
\end{Cor*}

\noindent
In fact, the Theorem might also be deduced from the statement of the Corollary, by using the six term exact sequence. 

Moreover, the Theorem can be used to prove Bott periodicity and thus the existence of the six-term exact sequence. This is the approach taken by Cuntz \cite{cuntz-bott}. (Cuntz defines $A_{\reals_+}$ algebraically, and his elegant deduction of the $K$-triviality of this C$^*$-algebra is quite different from the one we shall propose below.)  

\medskip\noindent
From the point of view of analysis and index theory, it seems natural to consider the multivariate generalisation of the Wiener--Hopf algebra and to study its $K$-theory. 

Thus, let $\Omega\subset\reals^n$ be a closed convex cone, which we assume to be \emph{pointed}, i.e.~$\Omega$ contains no affine line, and \emph{solid}, i.e.~$\Omega$ generates $\reals^n$ as a vector space. Then \emph{Wiener--Hopf operators} shall be the bounded operators on $\Lp[^2]0\Omega$ given by 
\[
(W_fg)(x)=\int_\Omega f(x-y)g(y)\,dy\mathfa g\in\Lp[^2]0\Omega\,,\,f\in\Lp[^1]0{\reals^n}\,,\,x\in\Omega\ .
\]
The C$^*$-algebra $A_\Omega$ of bounded operators generated by the $W_f$ will be called the \emph{Wiener--Hopf algebra}. This C$^*$-algebra and its relatives are the object of study of quite an extensive literature, and we refer the interested reader to the introduction of our joint paper with Troels Johansen \cite{alldridge-johansen1}, for a partial overview.

Just as in the $n=1$ case, there is an obvious symbol map $\sigma:A_\Omega\to\Ct[_0]0{\reals^n}\,$, which continues to be a surjective $*$-morphism for $n>1\,$. However, it is not to be expected that the kernel of $\sigma$ (the commutator ideal) equals the ideal of compact operators in this case. Rather, $A_\Omega$ has a composition series whose length is at most $n\,$. For the remainder of the paper, let us assume that $\Omega$ is \emph{polyhedral}, i.e.~finitely generated as a convex cone. Then the length of the composition series is exactly $n\,$, cf.~\cite{alldridge-johansen1}.

However, on the level of $K$- and even $KK$-theory, this distinction is invisible. Indeed, we shall prove in this paper the following theorem.

\begin{Th}
	Let $\Omega$ be a polyhedral cone. Then $A_\Omega$ is $KK$-contractible, and $A_\Omega/\knums$ and $S$ are $KK$-equivalent. 
\end{Th}

This Theorem was previously known only for a particular class of polyhedral cones called \emph{exhaustible}, and is due to Buyukliev \cite{buyukliev} in this case. He exploits the particular combinatorial structure of these cones to prove the Theorem via Mayer--Vietoris sequences and the exact six-term sequence. 

Arguably, in the general polyhedral case, the proof of $KK$-contractibility must take the whole combinatorial structure of an arbitrary polyhedral cone into account. In fact, the following result comes about as spin-off of our proof the above Theorem.

\begin{Th}
	Let $\Omega$ be a cone with polyhedral base $P\,$. Then the isomorphism class of $A_\Omega$ completely determines the combinatorial type of $P\,$, i.e., the lattice isomorphism class of its lattice of faces. 
\end{Th}

This is turn relies on the fact that the cellular differential of $P\,$, considered as a CW complex by considering each $j$-face as a $j$-cell, may be identified with the $d^1$ differential of the $E^1$ term of the Atiyah--Hirzebruch spectral sequence induced by the composition series alluded to above, a result which may be interesting in itself. We will describe this in detail below. 

\section{STRUCTURE OF THE WIENER--HOPF ALGEBRA}

\noindent
In this section, we review some results on the composition series of the Wiener--Hopf algebra $A_\Omega\,$, in particular, the construction and computation of certain `index maps'. These results are actually valid far beyond the case of polyhedral cones. However, restricting to polyhedral cones simplifies matters considerably, so we state them in this case only. The interested reader is referred to our joint papers with Troels Johansen \cite{alldridge-johansen1, alldridge-johansen2}, for the general case.

\subsection{Composition Series and Analytical Index Formula} 

Let $\Omega\subset\reals^n$ be a pointed and solid polyhedral cone. $\Omega$ is spanned by its exposed rays, and one may choose a set $E$ of generators of exposed rays contained in an affine hyperplane $H\,$. There exists a linear automorphism $L$ of $\reals^n$ such that $L(H)=1\times\reals^{n-1}\,$. Let $P$ be the convex hull of all $x$ such that $(1,x)\in L(E)\,$. Then $P$ is a convex polyhedron in $\reals^{n-1}\,$, and 
\[
L(\Omega)=\reals_+\cdot(1\times P)=\Set1{(\lambda,\lambda\cdot x)}{\lambda\sge0\,,\,x\in P\subset\reals^{n-1}}\ .
\]

Henceforth, we will omit reference to the linear automorphism $L\,$, and assume that $\Omega=\reals_+\cdot(1\times P)$ where the set $P\subset\reals^{n-1}$ is a convex polyhedron. This assumption is no loss of generality, since the C$^*$-algebras $A_\Omega$ and $A_{L(\Omega)}$ are isomorphic. 

For $j=-1,\dotsc,n-1\,$, let $f_j$ be the number of $j$-dimensional convex faces of $P$ (where the empty set is considered as the unique face of dimension $-1$). This somewhat annoying index shift is an artefact introduced by considering $j$-faces of $P$ as $(j+1)$-faces of $\Omega\,$, and will continue to trouble us in the following. 

\begin{Th*}[Muhly--Renault \cite{muhly-renault}]\label{th:whfiltration}
	There exists a finite filtration of $A_\Omega$ by ideals $I_0=0\subset I_1=\knums\subset\dotsm\subset I_{n+1}=A_\Omega$ \scth $I_{j+1}/I_j$ is a liminary C$^*$-algebra with spectrum $\{1,\dotsc,f_{j-1}\}\times \reals^j\,$. 
\end{Th*}

Both the ideals $I_j$ and the isomorphism of the subquotients $I_{j+1}/I_j$ with the algebras $\Ct[_0]0{\{1,\dotsc,f_{j-1}\}\times\reals^j}\otimes\knums\,$, $j<n\,$, and $\Ct[_0]0{\reals^n}\,$, $j=n\,$, respectively, are given quite explicitly, but we shall not need the precise formulae. 

Let $\partial_j:K^i_c(\{1,\dotsc,f_{j-1}\}\times\reals^j)\to K^{i+1}_c(\{1,\dotsc,f_{j-2}\}\times\reals^{j-1})\,$, $j=1,\dotsc,n\;$, be the $K$-theory connecting maps induced by the exact sequences
\[
\xymatrix{0\ar[r]&I_j/I_{j-1}\ar[r]&I_{j+1}/I_{j-1}\ar[r]^-{\sigma_j}&I_{j+1}/I_j\ar[r]&0}\ .
\]

Let $F_j$ be the set of $j$-dimensional faces of $P\,$, and let $\Omega_F\,$, for $F\in P\,$, be the face of $\Omega$ spanned by $F\,$. For any $A\subset\reals^n\,$, let $\Span0A$ denote the linear span. We may identify $\{1,\dotsc,f_{j-1}\}\times\reals^j$ with the trivial rank $j$ vector bundle 
\[
	\Sigma_j=\Set1{(F,y)\in F_{j-1}\times\reals^n}{y\in\Span0{\Omega_F}}
\]
over the finite base $F_{j-1}\,$. 

For any subset $A\subset\reals^n\,$, let $A^*=\Set0{y\in\reals^n}{\Scp0yA\sge0}$ be the dual cone, and let $A^\perp=\Set0{y\in\reals^n}{\Scp0yA=0}$ be the orthogonal complement of the linear span. Then define, for any face $F$ of $P\,$, 
\[
	\Omega_F^\circledast=\Set1{x\in\Span0{\Omega_F^\perp\cap\Omega^*}}{\Scp0xy\sge0\smathfa y\in \Omega_F^\perp\cap\Omega^*}\ .
\]
(This notation differs from \cite{alldridge-johansen1, alldridge-johansen2}.)

The continuous field of Hilbert spaces $(\Lp[^2]0{\Omega_F^\circledast})_{(F,y)\in\Sigma_j}$ naturally defines a Hilbert $\Ct[_0]0{\Sigma_j}$-module $\mathcal E_j\,$. The $*$-morphism $\sigma_j$ extends to a representation of $A_\Omega$ by adjointable endomorphisms of this Hilbert module. By these means, the map $\partial_j$ lends itself to an analytical expression, as follows. 

\begin{Th*}[A.--Johansen \cite{alldridge-johansen1}]
	Let $a\in M_N(I_{j+1}^+)$ represent the $K$-theory class $[\sigma_j(a)]\in K^1_c(\Sigma_j)\,$. Then $\sigma_{j-1}(a)$ is a Fredholm operator on the Hilbert $\Ct[_0]0{\Sigma_{j-1}}$-module $\mathcal E_{j-1}^N\,$, and 
	\[
		\partial_j[\sigma_j(a)]=\Index\sigma_{j-1}(a)\in K^0_c(\Sigma_{j-1})\ .
	\]
\end{Th*}

\subsection{$KK$-Theoretical Index Formula}

The finite set 
\[
	\mathcal P_j=\Set1{(E,F)\in F_{j-2}\times F_{j-1}}{E\subset F}
\]
may be considered as bibundle w.r.t.~the obvious projections $\xi:\mathcal P_j\to F_{j-2}$ and $\eta:\mathcal P_j\to F_{j-1}$ ($j\sge1$). The map
\[
	\eta^*\Sigma_j\to\xi^*\Sigma_{j-1}:(E,F,y)\mapsto(E,F,p_E(y))
\]
realises $\eta^*\Sigma_j$ as the trivial line bundle over the base $\xi^*\Sigma_{j-1}\,$. (Here, for $A\subset\reals^n\,$, $p_A$ denotes the orthogonal projection onto $\Span0A\,$.) Indeed, a nowhere vanishing section is given by the map 
\[
s_j:\xi^*\Sigma_{j-1}\to\eta^*\Sigma_j:(E,F,u)\mapsto(E,F,u+e_F(E))
\]
where the unit vector $e_F(E)\in\Omega_E^\perp\cap\Span0{\Omega_F}$ is given as follows. If $\check \Omega_F=\Omega_F^\perp\cap\Omega^*$ denotes the dual face of $\Omega_F\,$, then $\check\Omega_F^\perp\cap\Omega_E^\circledast$ is an extreme ray of $\Omega_E^\circledast\,$, and $e_F(E)$ is the unique unit vector contained in this ray.

The trivial line bundle $\eta^*\Sigma_j\to\xi^*\Sigma_{j-1}$ induces an isomorphism of $K$-groups $K_c^i(\eta^*\Sigma_j)\to K_c^{i+1}(\xi^*\Sigma_{j-1})\,$. It is given by multiplication by an invertible $KK$-theory element $y_j\in KK^1(\eta^*\Sigma_j,\xi^*\Sigma_{j-1})$ where for locally compact Hausdorff spaces $X$ and $Y\,$,  we write $KK^q(X,Y)=KK(\Ct[_0]0X,\Ct[_0]0{\reals^q\times Y})\,$. Another way to think about $y_j$ is that it is `fibre integration', i.e.~the inverse of the Thom isomorphism for the above line bundle. This depends on the choice of an orientation; we will go into detail further below. 

The only other ingredients needed for our index formula (at least in the polyhedral case) are the projections $p_\xi:\xi^*\Sigma_{j-1}\to\Sigma_{j-1}$ and $p_\eta:\eta^*\Sigma_j\to\Sigma_j$ induced, respectively, by $\xi$ and $\eta\,$. Since $\xi$ and $\eta$ have finite domain, $p_\xi$ and $p_\eta$ are proper, and thus induce $*$-morphisms $\pmb\xi:\Ct[_0]0{\Sigma_{j-1}}\to\Ct[_0]0{\xi^*\Sigma_{j-1}}$ and $\pmb\eta:\Ct[_0]0{\Sigma_j}\to\Ct[_0]0{\eta^*\Sigma_j}\,$, respectively. 

\begin{Th*}[A.--Johansen \cite{alldridge-johansen2}]
	Let $1\sle j\sle n\,$. As elements of $KK^1(\Sigma_{j-1},\Sigma_j)\,$, 
	\[
		\partial_j=\pmb\xi_*\pmb\eta^*y_j\ .
	\]
\end{Th*}

\section{COHOMOLOGICAL INDEX AND CELLULAR DIFFERENTIAL}

\subsection{Cohomological expression of the index}

If $X$ is a locally compact space, then there is a natural ring morphism $\ch:K_c^*(X)\to \bigoplus_{k=0}^\infty H^{*+2k}(X,\rats)\,$, called the \emph{Chern character}, which is rationally an isomorphism. 

Let $\pi:V\to X$ be an oriented real vector bundle. We have Thom isomorphisms
\[
	\vphi_V:K_c^*(X)\to K_c^{*+\rk V}(V)\nd\psi_V:H_c^*(X)\to H_c^{*+\rk V}(V)\ .
\]
In the special case that $V$ is \emph{trivial}, these are related via the Chern character, i.e.~$\ch\circ\vphi_V=\psi_V\circ\ch\,$. (In general, of course, the interplay is more subtle.) 
As is customary, we denote integration along the fibres of $\pi\,$, which is the inverse map $\psi_V^{-1}:H_c^{*+\rk V}(V)\to H_c^*(X)\,$, by $\pi_*\,$. 

In particular, applying this to $y_j^{-1}:K_c^*(\xi^*\Sigma_{j-1})\to K_c^{*+1}(\eta^*\Sigma_j)\,$, we immediately obtain the following cohomological expression of the index map $\partial_j\,$. 

\begin{Prop}
	For all $u\in K_c^i(\Sigma_j)\,$, we have
	\[
		\ch(\partial_j(u))=p_{\eta*}\pi_*p_\xi^*\ch(u)\mathtxt{in}\bigoplus_{k=0}^\infty H_c^{2k+i+1}(\Sigma_{j-1},\rats)
	\]
	where $\pi$ denotes the projection of the (trivial) line bundle $\eta^*\Sigma_j\to\xi^*\Sigma_{j-1}\,$. 
\end{Prop}

We remark that 
\[
K_c^i(\Sigma_j)=K_c^i(F_{j-1}\times\reals^j)=\begin{cases}0&i+j\equiv1\pmod 2\ ,\\\ints^{f_{j-1}}&i+j\equiv0\pmod2\ ,\end{cases}
\]
so that $K_c^*(\Sigma_j)$ has no torsion. In particular, $\ch$ is injective on $K_c^*(\Sigma_j)\,$. It is known that its image is the integral cohomology $H^*_c(\Sigma_j,\ints)\,$, cf.~\cite[Proposition 4.3]{hatcher-vbktheory}. In particular, $\partial_j$ is completely determined by its cohomological expression, and the latter is integral.

In we take the trivialisation of $\eta^*\Sigma_j\to\xi^*\Sigma_{j-1}$ to be given by the non-vanishing section $s_j$ defined above, then the orientation of this bundle is induced by choices of orientations of $\Sigma_{j-1}$ and $\Sigma_j\,$. These will be induced by choices of trivialisations of these bundles. (The triviality of the latter bundles is particular to the polyhedral situation.) As we shall presently see, the detailed inspection of these choices leads directly to the explicit expression of $\partial_j$ as a cellular differential.

\subsection{The Wiener--Hopf Index as a Cellular Differential}

The $n$-dimen\-sion\-al convex polytope $P$ can be considered as a finite CW complex by taking the $j$-faces to be the $j$-cells. (Topologically, $P$ is of course an $n$-cell, so there are simpler ways to consider  it as a CW-complex. However, our point of view captures the combinatorics of the face lattice.) The cellular complex is then $(H^0(F_j),d_j)$ where $H^0(F_j)$ is the free Abelian group generated by the $j$-faces. 

The vector bundle $\Sigma_j\to P_j$ is trivial, hence orientable, and we have a Thom isomorphism  $\psi_j:H^0(F_j)\to H^j_c(\Sigma_{j+1})$ given by the choice of an orientation. Let us make this choice explicit. A trivialisation of $\Sigma_j$ is given by the map 
\[
F_{j-1}\times\reals^j\to\Sigma_j:(F,y)\mapsto(F,A_Fy)
\]
where for each $F\in F_{j-1}\,$, $A_F:\reals^j\to\Span0{\Omega_F}$ is a linear isomorphism. An orientation of (the fibres of) $\Sigma_j$ is given by pulling back the standard orientation $\sigma^+=\eps\circ\det$ of $\reals^j$ to $\Span0{\Omega_F}$ along $A_F^{-1}$ to an orientation $\sigma_F\,$. The Thom isomorphism $\psi_j$ is given by the cup product with the Thom class $c_j\in H^j_c(\Sigma_j)$ which is determined by the condition 
\[
\int_{(\Span0{\Omega_F},\sigma_F)}c_j(F,\cdot)=1\mathfa F\in F_{j-1}\ .
\]

  The line bundle $\pi:\eta^*\Sigma_j\to\xi^*\Sigma_{j-1}$ is oriented by the choice of the non-vanishing section $s_j\,$. Observe that for each $(E,F)\in\mathcal P_j\,$, there exists a unique orientation $\sigma$ on $\Span0{\Omega_F}$ such that 
  \[
  \sigma(e_F(E),v_1,\dotsc,v_{j-1})=\sigma_E(v_1,\dotsc,v_{j-1})\mathfa v_1,\dotsc,v_{j-1}\in\Span0{\Omega_E}\ .
  \]
  We denote by $[E:F]=\pm1$ the unique sign such that $\sigma_F=[E:F]\cdot\sigma_E\,$. If $E\not\subset F\,$, we define $[E:F]=0\,$. 

\begin{Prop}
  For $0\sle j\sle d\,$, let $\tilde d_j:H^0(F_j)\to H^0(F_{j-1})$ be the map induced by the index map $\partial_{j+1}:K_c^{j+1}(\Sigma_{j+1})\to K_c^j(\Sigma_j)$ and the Thom isomorphisms $\psi_{j+1}\,$, $\psi_j\,$, via the relation 
  \[
  		\psi_j\circ\tilde d_j\circ\psi_{j+1}^{-1}\circ\ch=\ch\circ\,\partial_{j+1}\ .
  \]
  Then $\tilde d_j$ is given on generators by the formula 
  \[
  \tilde d_j(F)=\sum_{E\subset F}[E:F]\cdot E\ .
  \]
\end{Prop}

\begin{proof}
Consider a form $c_{j+1}\in\Gamma_c(\Sigma_{j+1},\wedge^{j+1}T^*\Sigma_{j+1})$ representing the Thom class in $H^{j+1}_c(\Sigma_{j+1})\,$. Then $p_\eta^*c_{j+1}$ is represented by 
\[
\eta^*\Sigma_{j+1}\to\wedge^{j+1}T^*\eta^*\Sigma_{j+1}:(E,F,u)\mapsto c_{j+1}(F,u)\ .
\]
Because we have the decomposition $\Span0{\Omega_F}=\reals\cdot e_F(E)\oplus\Span0{\Omega_E}\,$, the Fubini theorem gives 
  \[
  \int_{(\Span0{\Omega_E},\sigma_E)}\int_\reals p_\eta^*c_{j+1}(E,F,\cdot+te)(e,\cdot)\,dt=[E:F]\cdot\int_{(\Span0{\Omega_F},\sigma_F)}c_{j+1}(F,\cdot)=[E:F]\ ,
  \]
  where we write $e=e_F(E)\,$. Since this condition characterises the Thom class $c_j$ up to the factor $[E,F]\,$, 
  \[
  \pi_*p_\eta^*c_{j+1}(E,F,\cdot)=[E:F]\cdot c_j(E,\cdot)\mathtxt{in}H_c^j(\Span0{\Omega_E})\ .
  \] 
  We find 
  \[
  p_{\xi*}\pi_*p_\eta^*c_{j+1}(E,u)=\sum_{F\subset E}\pi_*p_\eta^*c_{j+1}(E,F,u)=\sum_{F\subset E}[E:F]\cdot c_j(E,u)\ .
  \]
  Applying the cup product, the conclusion follows. 
\end{proof}

In order to see that the maps $\tilde d_j$ coincide with the cellular differentials $d_j\,$, let us explicitly describe $d_j\,$. To that end, we construct attaching maps. First, for any convex set $C\subset\reals^j$ containing the origin in its interior, let $\mu_C$ be its Minkowski gauge. If, more generally, the convex set $C$ has non-void interior and $b_C$ is its barycentre, then the map
\[
\vphi_C:\reals^j\to\reals^j\ , \ 
\vphi_C(x)=\begin{cases}\frac{\mu_{C-b_C}(x-b_C)}{\Norm0{x-b_C}}\cdot(x-b_C)&x\neq b_C\ ,\\
0&x=b_C\ ,\end{cases}
\]
is a homeomorphism inducing homeomorphisms $C\to\ball^j$ and $\partial C\to\sph^{j-1}$ where the latter has degree $1\,$. 

Next, for any $j$-face $F\in F_j\,$, let the linear isomorphism $A_F:\reals^{j+1}\to\Span0{\Omega_F}$ used above to trivialise $\Sigma_{j+1}$ be chosen such that $A_F(0\times\reals^j)$ is the linear subspace parallel to the affine span of $F\,$, and \scth $\tilde F:=A_F^{-1}(F-b_F)\subset\ball^j\,$. We then define the attaching map for the $j$-cell associated with $F$ by
\[
	\phi_F:\ball^j\to F\subset P:x\mapsto A_{\vphantom{\tilde F}F}^{\vphantom{-1}}\vphi_{\tilde F}^{-1}(x)+b_F
\]

If for any pair $(E,F)\in F_{j-1}\times F_j$ the number $c_{EF}\in\ints$ denotes the degree of the composite
 \[
  \xymatrix{%
   \sph^{j-1}\ar[r]^-{\phi_F}&\partial F\ar[r]&P/(P\setminus E^\circ)=E/\partial E\ar[r]^-{\phi_E^{-1}}&\ball^{j-1}/\sph^{j-2}\ar[r]&\sph^{j-1}}
  \]
  where the rightmost map is the standard one ($x\mapsto(ux,2\Norm0x-1)$ where $u$ is suitably chosen), then the cellular differential $d_j:H^0(F_j)\to H^0(F_{j-1})$ is defined on generators by
  \[
  d_j(F)=\sum_{E\subset F}c_{EF}\cdot E\ .
  \]

\begin{Prop}[degree-coeff]
  We have $c_{EF}=[E,F]$ for any pair $(E,F)\in F_{j-1}\times F_j\,$. In particular, the cellular differential $d_j$ coincides with $\tilde d_j\,$. 
\end{Prop}

\begin{proof}
Let $H_0:\partial F\to\sph^{j-1}$ be the composite
\[
  \xymatrix{%
   \partial F\ar[r]&P/(P\setminus E^\circ)=E/\partial E\ar[r]^-{\phi_E^{-1}}&\ball^{j-1}/\sph^{j-2}\ar[r]&\sph^{j-1}}\ .
  \]
  Next, let $1\sge t>0\,$. Let $e$ be a positive multiple of the projection of $e_F(E)$ onto the subspace $0\times\reals^n=\Span0P\,$, and define 
  \[
  F_s=\Set1{x\in F}{\Rscp0{x-b_E}e\sle s}\mathfa s\in[0,1]\ ,
  \]
  the elements of `height' $\sle s$ over $E\,$. By definition of $e\,$, $F_0=E\,$. We may take $e$ to be normalised in such a way that $F_1=F\,$. 
  
  We wish to define a map $H_t:\partial F_t\to\sph^{j-1}\,$. To that end, 
  define an affine map 
  \[
  A_t:\reals\times\reals^{j-1}\to\reals^n\mathtxt{by}A_t(2s-t,x)=s\cdot e+A_Ex+b_E\ .
  \]
  Then $\tilde F_t=A_t^{-1}(F_t)$ is a compact convex subset of $\reals^j$ containing $0$ in its interior. Define 
  \[
  H_t:F_t\cap\partial F\to\sph^{j-1}\mathtxt{by}H_t=f_t\circ\vphi_t^{-1}\circ A_t^{-1}\ ,
  \]
  where $\vphi_t=\vphi_{\tilde F_t}$ and for $r\in[0,1]\,$, $f_r:\sph^{j-1}\to\sph^{j-1}$ is given by 
  \[
  f_r(s,x)=\Parens1{ux,\min\Parens1{1,-1+2\tfrac{s+1}{r+1}}}\ ,
  \]
  for suitably chosen $u\,$. The map $f_r$ maps all points of the sphere of height $\sge r$ to the `north pole' $e_j=(0,\dotsc,0,1)\,$. 
  
  The set 
  $B=\partial F_t\cap\Set1x{\Rscp0{x-b_E}e=t}$ bounds a flat in $F_t\,$. Since $\tilde F_t\subset\ball^j\,$, the elements of $\vphi_t^{-1}(A_t^{-1}(B))$ have 
  $j$th coordinate $\sge t\,$. Thus, $H_t$ maps $B$ to $e_j\,$, and hence extends to all of $\partial F$ by sending $\partial F\setminus F_t$ to $e_j\,$. Moreover, 
  $H_t\,$, together with $H_0\,$, form a homotopy. We conclude
  \[
  c_{EF}=\deg H_0\circ\phi_F=\deg H_1\circ\phi_F=\mathop{\mathrm{sign}}\det((e,A_E)^{-1}A_F)=[E:F]\ ,
  \]
  which proves our claim.
\end{proof}

\section{PROOF THE MAIN THEOREM}

\subsection{Invariance of the Combinatorial Type of $P$}

Recall that the \emph{combinatorial type} of the convex polyhedron $P$ is the lattice isomorphism class of the lattice of convex faces of $P\,$. The \emph{$f$-vector} of $P$ is the vector $(f_0,\dotsc,f_n)$ whose component $f_j$ is the numbers of $j$-faces. 

The following theorem is a somewhat surprising if simple consequence of \thmref{Prop}{degree-coeff}. 

\begin{Th}
  Let $\Omega$ be a convex cone with polyhedral base $P\,$. Then the isomorphism class of $A_\Omega$ determines the combinatorial type of $P\,$. I.e., if $\Omega'$ is another cone with polyhedral base $P'\,$, and $A_\Omega$ and $A_{\Omega'}$ are isomorphic, then $P$ and $P'$ have isomorphic face lattices. 
\end{Th}

\begin{proof}
  The ideals in the filtration $(I_j)$ of $A_\Omega$ from \thmref{Th}{whfiltration} are recursively characterised by the property that $I_{j+1}/I_j$ is the largest liminary ideal of $A_\Omega/I_j$ with Hausdorff spectrum. Thus, the $f$-vector of $P$ and the index maps, and thus the maps $d_j\,$, are uniquely determined up to a choice of orientations. In particular, the absolute values $\Abs0{[E:F]}$ are uniquely determined for any pair of faces $(E,F)$ where $E$ is of codimension one in $F\,$. But these numbers determine the lattice order. Hence the assertion. 
\end{proof}

\subsection{The $KK$-contractibility of $A_\Omega$}

As C$^*$-algebras with finite ideal filtrations with subquotients stably isomorphic to multiples of $\Ct[_0]0{\reals^n}\,$, $A_\Omega$ and $A_\Omega/\knums$ belong to the bootstrap category $\mathcal N$ and therefore obey the UCT \cite[Definition 22.3.4, Theorem 23.1.1]{blackadar-ktheory}. In order to determine their $KK$ equivalence class, it suffices to compute their $K$-theory. 

The C$^*$-algebra $A_\Omega$ has a filtration by ideals $(I_j)\,$, and since $I_1=\knums\,$, $A_\Omega/\knums$ has the filtration by ideals given by $(I_j/I_1)\,$. For any C$^*$-algebra $A$ filtered by ideals $(I_j)\,$, Schochet has introduced an Atiyah--Hirzebruch type homology spectral spectral sequence $(E^r_{p,q})$ which converges to the $K$-theory of $A\,$. In the case of $A=A_\Omega\,$, by \cite[Theorem 2.1]{schochet-topmeth1}, its $E^1$ term is 
  \[
  E^1_{p,q}=K_{p+q}(I_p/I_{p-1})=
  \begin{cases}
    H^0(F_{p-2})&q\equiv1\pmod2\ ,\\
    0&q\equiv0\pmod2\ .
  \end{cases}
  \]
Hence, every other column of $E^1$ is zero. Hence, $E^r$ abuts to $E^2\,$. 

\begin{Prop}
  The homology spectral sequence $E_{p,q}^r$ in $K$-theory induced by the filtration $(I_j)$ abuts to its $E^2$ term, which is zero. In particular, $K_*(A_\Omega)=0$ and $K_*(A_\Omega/\knums)=K_*(\cplxs)\,$. 
\end{Prop}

\begin{proof}
By definition, the $d^1$ differential is the composite
\[
\xymatrix{%
E^1_{p,1}=K_{p+1}(I_p/I_{p-1})\ar[r]^-{\partial}&K_p(I_{p-1})\ar[r]&K_p(I_{p-1}/I_{p-2})=E^1_{p-1,1}}
\]
where the first map is the boundary map in the exact six-term sequence in $K$-theory, and the second is induced by the quotient map  $I_{p-1}\to I_{p-1}/I_{p-2}\,$. Considering 
the commutative diagram with exact rows, 
\[
\xymatrix{%
0\ar[r]&I_{p-1}\ar[r]\ar[d]&I_p\ar[r]\ar[d]&I_p/I_{p-1}\ar@{=}[d]\ar[r]&0\\
0\ar[r]&I_{p-1}/I_{p-2}\ar[r]&I_p/I_{p-2}\ar[r]&I_p/I_{p-1}\ar[r]&0
}
\]
it follows from the naturality of connecting maps that $d^1$ is also given by the connecting map for the lower line. This is just the map $\partial_{p-1}\,$. We have already noted that $\partial_{p-1}$ is uniquely determined by its cohomological expression, and the latter gives the cellular differential $d_{p-2}\,$. Thus, $(E^1_{p,1},d^1_p)$ 
is up to a shift, just the augmented cellular chain complex of the CW complex $P\,$. Since $P$ is contractible, this complex is exact. Hence, $E^2=0\,$, and the first statement follows.

To complete the proof, observe that dividing by $I_1=\knums$ corresponds to removing the augmentation from the cellular complex. The resulting complex has cohomology concentrated in degree zero, and $H^0(F_0)=\ints\,$. (Alternatively, use the exact six-term sequence.)
\end{proof}

\begin{Cor}
	The C$^*$-algebra $A_\Omega$ is $KK$-contractible, and $A_\Omega/\knums$ and $S$ are $KK$-equivalent. 
\end{Cor}

\begin{Cor}
	The isomorphism $K_1(A_\Omega/\knums)\to\ints$ given by computing the numerical index of Fredholm Wiener--Hopf operators, is an isomorphism.
\end{Cor}

\begin{proof}
The groups $K_1(A_\Omega/\knums)$ and $\ints$ are isomorphic, and Buyukliev \cite{buyukliev} has constructed a Fredholm Wiener--Hopf operator of index one. 
\end{proof}

\medskip\emph{Acknowledgements.} 
Part of this work was conducted while the author was a visitor at the Institut Henri Poincar\'e, Paris. The author wishes to thank the institute for its hospitality and the organisers of the Special Trimester on `Groupoids, Stacks in Physics and Geometry' for their support.


\begin{thebibliography}{99}

\bibitem{alldridge-johansen1} \textsc{A.~Alldridge, T.R.~Johansen}, Spectrum and Analytical
	Indices of the C$^*$-Algebra of Wiener--Hopf Operators, \textit{J. Funct. Anal.} 
	\textbf{249}(2007), no.~2, 425--453.

\bibitem{alldridge-johansen2} \textsc{A.~Alldridge, T.R.~Johansen}, An Index Theorem for Wiener--Hopf Operators, \textit{Adv.~Math.} \textbf{218}(2008), no.~1, 163--201.

\bibitem{blackadar-ktheory} \textsc{B.~Blackadar}, \textit{$K$-Theory for Operator Algebras}, 2nd edition, MSRI Publications \textbf{5}, Cambridge University Press, Cambridge, 1998.
	
\bibitem{buyukliev} \textsc{N.P.~Buyukliev}, \textit{$K$-Theory of the C$^*$-Algebra of Multivariable Wiener--Hopf Operators Associated with some Polyhedral Cones in $\mathbb R^n$}, Annuaire Univ. Sofia Fac. Math. Inform. \textbf{91}(1997), no.~1--2, 115--125. 

\bibitem{cuntz-bott} \textsc{J.~Cuntz}, $K$-Theory and C$^*$-Algebras, pp.~55--79 in:~A.~Bak (ed.), \textit{Algebraic $K$-Theory, Number Theory, Geometry and Analysis (Bielefeld,
	1982)}, Lect. Notes Math. \textbf{1046}, Springer-Verlag, Berlin, 1984.

\bibitem{hatcher-vbktheory} \textsc{A.~Hatcher}, \textit{Vector Bundles and $K$-Theory}, Version 2.0, 2003.

\bibitem{muhly-renault} \textsc{P.S.~Muhly, J.N.~Renault}, C$\sp{\ast}$-Algebras of Multivariable 
	{W}iener--{H}opf Operators, \textit{Trans. Amer. Math. Soc.} \textbf{274}(1982),
	no.~1, 1--44.

\bibitem{schochet-topmeth1} \textsc{C.~Schochet}, Topological Methods for C*-Algebras I: Spectral Sequences, \textit{Pac. J. Math.}, \textbf{96}(1981), no.~1, 193--211.

\end{thebibliography}
\end{document}